\documentclass{amsart}
\usepackage[T1]{fontenc}
\usepackage{color}
\usepackage{textcomp}
\usepackage{subfig}
\usepackage{amsmath}
\usepackage{amssymb}
\usepackage{amsthm}
\usepackage[english]{babel}
\usepackage{graphicx}
\usepackage{tikz}
\usetikzlibrary{trees,intersections,calc,positioning,arrows}
\usepackage{hyperref}

\theoremstyle{definition}
\newtheorem{definition}{Definition}
\newtheorem{theorem}{Theorem}
\newtheorem{lemm}[theorem]{Lemma}

\newtheorem{corollary}[theorem]{Corollary}

\newcommand{\LS}{\ensuremath{\underset{n=1}{\overset{\infty}{\cap}} \, {\underset{k=n}{\overset{\infty}{\cup}}}\,}}

\newcommand{\LSM}{\ensuremath{\underset{n=1}{\overset{\infty}{\cap}} \, {\underset{m=n}{\overset{\infty}{\cup}}}\,}}

\newcommand{\ceilup}[1]{\left\lceil #1 \right\rceil}
\newcommand{\floordown}[1]{\left\lfloor #1 \right\rfloor}

\begin{document}
\title[Hausdorff Dimension of Measures for Non-Uniquely Ergodic IETs]{On the Hausdorff Dimension of Measures for a Non-Uniquely Ergodic Family of Interval Exchange Transformations}
\author{Aleksei Kobzev}
\address{Faculty of Mathematics, National Research University Higher School of Economics, Vavilova St. 7, 112312 Moscow, Russia}
\email{akobzev@hse.ru}

\begin{abstract}
In this paper, based on a construction by J. Fickenscher, we construct a family of non-uniquely ergodic interval exchange transformations on $n$ intervals with the maximal possible number of measures, $\floordown{\frac{n}{2}}$. Subsequently, we generalize J. Chaika's result on estimating the Hausdorff dimension of the two measures from M. Keane's example to the family of interval exchange transformations with the maximal possible number of measures.
\end{abstract}
\maketitle

\section{Introduction}
\subsection{Historical Overview}

Interval exchange transformations arise as first return maps to a transversal in the case of an orientable measured foliation on an orientable surface or, for example, as the first return map to a transversal for a billiard flow in a rational polygon.

Interval exchange transformations were first mentioned by V. Arnold in 1963 in his work \cite{Arn}, where he considered IETs on three intervals. The general case was first described by V. Oseledets in 1968 in \cite{Odc}.

A transformation is called minimal if the orbit of every point is dense. In his 1975 paper \cite{Kn75}, M. Keane formulated a sufficient condition for minimality: an interval exchange transformation satisfies the Keane condition if it never happens that applying the transformation to one of the boundary points maps it to another boundary point. In the same paper, M. Keane proved that almost all interval exchange transformations are minimal and conjectured, which was later independently proven by H. Masur and W. Veech in \cite{M82} and \cite{V82} respectively, that almost all interval exchange transformations are uniquely ergodic.

The next step in the study of interval exchange transformations was the construction of non-uniquely ergodic examples. Examples of minimal interval exchange transformations admitting two measures were first constructed by Keynes - Newton \cite{KN} and M. Keane \cite{Kn77}. M. Keane's construction will be used later in our work. The method is based on proving the existence of orbits with asymptotically different distributions, which leads to two distinct ergodic measures $\mu$ and $\nu$.

Subsequently, an upper bound on the number of invariant ergodic measures was obtained by A. Katok \cite{Kt}: an interval exchange transformation cannot have more than $m$ invariant ergodic measures, where $m$ is the genus of the surface constructed from the IET, i.e., $\lfloor \frac{n}{2} \rfloor$, where $n$ is the number of intervals. The first proof of the effectiveness of this bound was obtained by E. Sataev \cite{St}. This result was also obtained by J. Fickenscher \cite{Fk}.

One of the main tools for studying interval exchange transformations is Rauzy induction. To construct non-uniquely ergodic examples, one can use cycles on the Rauzy diagram. J. Fickenscher and J.-C. Yoccoz, in their works \cite{Fk} and \cite{Yc} respectively, described constructions of cycles on the Rauzy diagram that allow for the construction of examples of interval exchange transformations possessing two or more invariant ergodic measures. The example presented by M. Keane was studied in more detail by J. Chaika. In \cite{Ch8} and \cite{Ch11}, estimates for the Hausdorff dimension of the invariant measures from M. Keane's example were obtained.

\subsection{Main Result of the Work}

In this paper, our analysis centers on a significant family of interval exchange transformations constructed by J. Fickenscher in \cite{Fk}. These transformations, built using a specific closed path $\gamma^{a,c}$ on the Rauzy diagram, are remarkable for possessing the maximal possible number of invariant measures. This path is used to construct a sequence of parameters $\{a_i, c_i\}$ defining an IET on $n$ intervals, where the integer parameters $a$ and $c$ correspond to loops on the diagram. These parameters allow us to compute the transition matrix $\Theta_{c,a}$ between the initial and final lengths obtained through a cyclic path on the Rauzy diagram. 

We begin by reformulating Fickenscher's result into a more precise statement, presented below as Theorem \ref{onto}, which is tailored for our analytical approach. Our contribution lies in a deeper quantitative analysis of Fickenscher's construction. We will study the dynamical properties of the transition matrix $\Theta_{c,a}$, which, in turn, depends on the parameters $c$ and $a$. Subsequently, we examine the invariant measures that arise as limits of sequences of distributions. For each $k \in \{1, \dots, d\}$, we define a sequence of vectors where the limiting expression will have the form:
$$\lim_{m \to +\infty} \overline{\Theta}_{c_1,a_1}\overline{\Theta}_{c_2,a_2}...\overline{\Theta}_{c_m,a_m} e_k.$$
The initial condition corresponds to the canonical basis vector $e_k$.

Let $\lambda_k$ denote the invariant measure, which is the aforementioned limit. The measure $\lambda_k$ represents the distribution associated with an initial state concentrated in the $k$-th interval. The following theorem will be proven in this work, which states that for $k \in \{3, 5, \dots, 2\floordown{\frac{n}{2}}-1,n\}$, the corresponding measures $\lambda_k$ are pairwise distinct, which is a sufficient condition for non-unique ergodicity.


\begin{theorem}\label{onto}
    Let an interval exchange transformation be defined by a sequence of parameters $\{a_i, c_i\}$ generated by the Rauzy cycle $\gamma^{a,c}$, such that $c_{i+1} = pa_i = p^2c_i$ for some integer $p \geq n + 1$. Such a transformation is non-uniquely ergodic.

    Specifically, it possesses $\floordown{\frac{n}{2}}$ distinct invariant ergodic measures, found as the limits
    $$\lambda_j = \lim_{m \to +\infty} \overline{\Theta}_{c_1,a_1}\overline{\Theta}_{c_2,a_2}...\overline{\Theta}_{c_m,a_m} e_j, \text{ for } j \in \{3, 5, \dots, 2\floordown{\frac{n}{2}}-1, n\}.$$
\end{theorem}
$$ $$

The obtained estimates allow us to study in detail the combinatorial properties of the interval exchange transformation from the family defined in Theorem \ref{onto}, which enables us to derive a pairwise estimate of the Hausdorff dimension for the constructed invariant measures. Using the approach described by J. Chaika, we prove the following theorem for the example with $n$ intervals and $\floordown{\frac{n}{2}}$ distinct invariant ergodic measures.

\begin{theorem}\label{onto2}
In the setting of Theorem \ref{onto}, \\ for any distinct $i,j \in \{3,5,\dots,2\floordown{\frac{n}{2}}-1, n\}$, the following estimate holds:
$$ 
\underset{k \to \infty}{\liminf} \log_{\lambda_j(I_i^{(k)})} (\lambda_i(I_i^{(k)})) \leq \dim_H(\lambda_i,d_{\lambda_j}) \leq \underset{k \to \infty}{\liminf} \log_{\lambda_j(I_i^{(k)})} b^{-1}_{k,i},
$$
where $b_{k,i}=|\Theta_1 \cdots \Theta_k e_i|.$
\end{theorem}

In particular, for the sequence $\{a_i, c_i\}$ ($c_{i+1} = pa_i = p^2c_i$, $p \geq n + 1$), the following holds:
\begin{corollary}\label{onto3}
In the setting of Theorem \ref{onto}, \\ for any distinct $i,j \in \{3,5,\cdots,2\floordown{\frac{n}{2}}-1,n\}$, the following estimate holds
$$\dim_H(\lambda_i,d_{\lambda_j}) = 1$$
\end{corollary}

\subsection{Organization of the paper}

In Section \ref{sec:2}, we provide the necessary background, including basic definitions and the key mathematical tools used throughout the paper.

Section \ref{sec:3} is dedicated to a detailed description of J. Fickenscher's construction of a family of non-uniquely ergodic IETs. We then derive new estimates for this family, which serve as the foundation for our subsequent analysis. These estimates allow us to investigate the combinatorial properties of the considered IETs in greater depth. Specifically, by applying these bounds and employing M. Keane's method from \cite{Kn75}, we provide the proof of Theorem \ref{onto}.

Finally, in Section \ref{sec:4}, we generalize the results of J. Chaika \cite{Ch11}. We establish additional estimates regarding the combinatorial structure of the IET family, which enable us to calculate the Hausdorff dimension of each invariant measure. These results complete the proofs of Theorem \ref{onto2} and Corollary \ref{onto3}

\subsection{Acknowledgments}
I am deeply grateful to A. Skripchenko for her mentorship and unwavering support during the preparation of this paper. This work is an output of a research project implemented as part of the Basic Research Program at the National Research University Higher School of Economics.

\subsection{Some notations}

\begin{description}
    \item[$\Theta$] The transition matrix according to the Rauzy graph. \\ By default, $\Theta = \Theta_{a,c}$. We also assume by default that $\Theta x^{j+1} = \Theta_{a^j,c^j} x^{j+1}$. The expression $(\Theta x^{j+1})_k$ denotes the $k$-th component of the vector obtained by applying the matrix $\Theta$ to the vector $x^{j+1}$.
    \item[$p, c, a$] Positive constants used in the estimates. In the proof, it is assumed that $p$ is a sufficiently large parameter for the inequalities to hold. The constants $a, c$ define the elements of the matrix $\Theta$. We will later define a sequence of parameters $\{a\}_i$ and $\{c\}_i$ which will be related by:
    $$
    c^{i+1} = p a^i = p^2 c^i
    $$
    \item[$x^j$] A vector with $n$ components at the $j$-th step, $x^j = (x^j_1, x^j_2, \dots, x^j_n)$: $$\sum_i x_i^{j+1} = 1.$$
    \item[$|\Theta x^{j+1}|$] This is the $L_1$-norm, i.e., $|\mathbf{v}| = \sum_i v_i$.
    \item[$\overline{\Theta}$] Denotes the following: $\overline{\Theta} = \frac{\Theta}{|\Theta|}$. \\ If no vector is specified, it defaults to $\frac{\Theta x^{j+1}}{|\Theta x^{j+1}|}$. The denominator is as in the previous item.
    \item[$\lambda_k(\cdot)$] Notation for measures. It will be shown in the paper that each such measure is an invariant ergodic probability measure and that measures with different indices are distinct. An important notational convention is:
    $$
    \lambda_k(I^k) = x^k, \ \lambda_k(I^k_i) = x^k_i.
    $$
    \item[$I_k, I^j_i$] Notations for subintervals. Conventionally, $I^j_i$ means the $i$-th subinterval of the interval at the $j$-th step of the first return map.
    \item[$\ceilup{\cdot}, \floordown{\cdot}$] Notations for ceiling and floor, respectively.

    \item[$b_{k, i}$] Time to the first return of $I_i^k$ to $I^k$. $b_{k, i} = |\Theta_1 ... \Theta_k e_i|_1$.

    \item[$O(I_i^k)$] Images of $I_i^k$ under iteration $T$ up to the first return to $I^k$: $\cup_{j=1}^{b_{k,j}} T^j(I_i^k)$.
\end{description}

\section{Interval Exchange Transformations}\label{sec:2}

In this section, we will present the formal definitions and tools that we will use to prove the results presented in this paper.

\subsection{Basic definitions}

\begin{definition} 
Let $X = [0, 1)$ and let $X$ be partitioned into $X_1, ... , X_n$, where $n$ is an integer greater than two. 
Let $\alpha = (\alpha_1, ..., \alpha_n)$ be a probability vector and $\tau$ be a permutation of the numbers $\{1, ..., n\}$.

Then we set $\beta_0 = 0$ and $\beta_i = \sum\limits_{j=1}^i \alpha_j$. Now let $X_i = [\beta_{i-1}, \beta_i)$. 

Then we define $\alpha^{\tau}$ and the corresponding $\beta^{\tau}_i$ and $X^{\tau}_i$:

$$\alpha^{\tau} = (\alpha_{\tau^{-1}(1)}, ..., \alpha_{\tau^{-1}(n)}), \ \beta^{\tau}_i = \sum\limits_{j=1}^i \alpha^{\tau}_j, \ X^{\tau}_i = [\beta^{\tau}_{i-1}, \beta^{\tau}_i).$$

Now we define the map $T: X \to X$ as follows:

$$T(x) = x - \beta_{i-1} + \beta^{\tau}_{\tau(i)-1}, \forall x \in X_i, 1 \leq i \leq n.$$

$T$ maps the interval $X_i$ bijectively to the interval $X^{\tau}_{\tau(i)}$.

The map $T$ is called an interval exchange transformation and is defined by the pair $(\alpha, \tau).$\\

$$\begin{tikzpicture}
\draw[ultra thick, red] (0,0) -- (1,0);
\draw[ultra thick, green] (1,0) -- (6,0);
\draw[ultra thick, blue] (6,0) -- (10,0);
\node at (0.5,0.5) {\color{black} $X_1$};
\node at (3.5,0.5) {\color{black} $X_2$};
\node at (8,0.5) {\color{black} $X_3$};
\draw[ultra thick, red] (9,-0.5) -- (10,-0.5);
\draw[ultra thick, green] (0,-0.5) -- (5,-0.5);
\draw[ultra thick, blue] (5,-0.5) -- (9,-0.5);
\node at (9.5,-1) {\color{black} $X_1$};
\node at (2.5,-1) {\color{black} $X_2$};
\node at (7,-1) {\color{black} $X_3$};
\end{tikzpicture}$$

\begin{center}
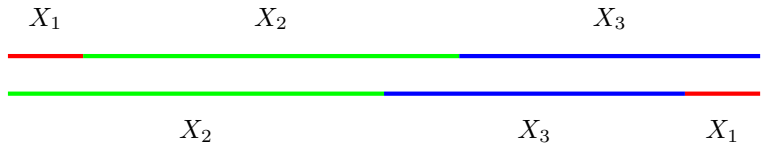
\textbf{Figure 1:} Interval exchange transformation on three subintervals.
\end{center}
\end{definition} 

\subsection{Rauzy Induction} 

\begin{definition}
Rauzy induction is the main renormalization procedure that allows us to construct a new IET from an initial one, which has the exact same orbit structure but is defined on a smaller support interval. Rauzy induction was first introduced by G. Rauzy for IETs in his 1979 paper \cite{Rz}.\\

Right Rauzy Induction:
Consider an interval exchange transformation $T: I \to I$ on $n$ subintervals of the interval $I = [0, 1)$. We take the rightmost subintervals in the domain and the range of the transformation. Then we induce $T$ on $[0, 1 - \lambda)$, where $\lambda = \min(\lambda_{a_n}, \lambda_{a_i})$ is the length of the shorter of the rightmost subintervals of $I$ and $TI$.\\
If $\lambda_{a_n} > \lambda_{a_i}$, the action of Rauzy induction is denoted as "0".\\
If $\lambda_{a_i} > \lambda_{a_n}$, the action of Rauzy induction is denoted as "1".\\

We will also say that the longer interval is the "winner" and the shorter one is the "loser".
\end{definition} 

The diagrams below show the transition from the permutation $\pi$ of the transformation $T$ to the permutation $\pi'$ of the transformation $T'$ under Rauzy induction of types "0" and "1". \\

\begin{center}

\begin{tikzpicture}[
    ->,
    >=stealth,
    shorten >=5pt,
    auto,
    node distance=7cm,
    semithick
]

\tikzset{
    state/.style={fill=white, draw=none, text=black, rounded corners, inner sep=4pt}
}

\node[state] (A) {$\begin{Bmatrix}
a_1&...&&&...&a_n\\
...&...&a_n&a_j&...&a_i
\end{Bmatrix}$};
\node[state] (B) [right of=A] {$\begin{Bmatrix}
a_1&...&&&...&a_n\\
...&a_n&a_i&a_j&...&...
\end{Bmatrix}$};

\path[->] (A) edge node {0} (B);

\end{tikzpicture}

\begin{tikzpicture}[
    ->,
    >=stealth,
    shorten >=5pt,
    auto,
    node distance=7cm,
    semithick
]

\tikzset{
    state/.style={fill=white, draw=none, text=black, rounded corners, inner sep=4pt}
}

\node[state] (A) {$\begin{Bmatrix}
a_1&...&a_i&a_{i+1}&...&a_n\\
...&...&&&...&a_i
\end{Bmatrix}$};
\node[state] (B) [right of=A] {$\begin{Bmatrix}
a_1&...&a_i&a_n&a_{i+1}&...\\
...&...&&&...&a_i
\end{Bmatrix}$};

\path[->] (A) edge node {1} (B);

\end{tikzpicture}

\end{center}

The implementation of right Rauzy induction looks as follows:

$$
\begin{tikzpicture}[scale=0.8]

\draw[ultra thick, red] (0,0) -- (1,0);
\draw[ultra thick, blue] (1,0) -- (3,0);
\draw[ultra thick, green] (3,0) -- (7,0);

\draw[ultra thick, red] (9,0) -- (13,0);
\draw[ultra thick, blue] (13,0) -- (15,0);
\draw[ultra thick, green] (15,0) -- (16,0);

\draw[ultra thick, green] (0,-0.6) -- (4,-0.6);
\draw[ultra thick, blue] (4,-0.6) -- (6,-0.6);
\draw[ultra thick, red] (6,-0.6) -- (7,-0.6);

\draw[ultra thick, green] (9,-0.6) -- (10,-0.6);
\draw[ultra thick, blue] (10,-0.6) -- (12,-0.6);
\draw[ultra thick, red] (12,-0.6) -- (16,-0.6);

\draw[ultra thick, red] (0,0-1.7) -- (1,0-1.7);
\draw[ultra thick, blue] (1,0-1.7) -- (3,0-1.7);
\draw[ultra thick, green] (3,0-1.7) -- (6,0-1.7);

\draw[ultra thick, red] (9,0-1.7) -- (12,0-1.7);
\draw[ultra thick, blue] (13,0-1.7) -- (15,0-1.7);
\draw[ultra thick, green] (12,0-1.7) -- (13,0-1.7);

\draw[ultra thick, green] (0,-0.6-1.7) -- (3,-0.6-1.7);
\draw[ultra thick, blue] (4,-0.6-1.7) -- (6,-0.6-1.7);
\draw[ultra thick, red] (3,-0.6-1.7) -- (4,-0.6-1.7);

\draw[ultra thick, green] (9,-0.6-1.7) -- (10,-0.6-1.7);
\draw[ultra thick, blue] (10,-0.6-1.7) -- (12,-0.6-1.7);
\draw[ultra thick, red] (12,-0.6-1.7) -- (15,-0.6-1.7);

\draw[dashed] (6,-3) -- (6,0.3);
\draw[dashed] (15,-3) -- (15,0.3);

\draw[ultra thick, ->, black] (3,-0.8) -- (3,-1.3);
\draw[ultra thick, ->, black] (12,-0.8) -- (12,-1.3);
\end{tikzpicture}
$$
\begin{center}
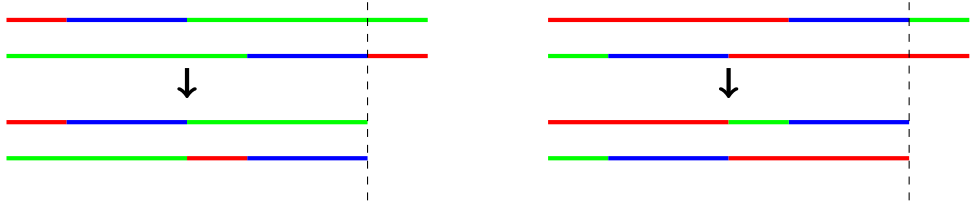
\textbf{Figure 2:} Two possible options of Rauzy induction.
\end{center}

\subsection{Hausdorff Dimension} 

\begin{definition}
Let $X = (X, d)$ be a metric space. Then for any $\alpha > 0$, $H_{\alpha}(*)$ denotes the Hausdorff $\alpha$-measure on $X$. $H_{\alpha}(*)$ is an outer measure defined for all subsets of $X$ as follows. First, for $A \subseteq X$ and $\epsilon > 0$, consider
$$
H_{\alpha, \epsilon}(A) = \inf \sum (r_i)^{\alpha},
$$
where the infimum is taken over all countable covers of $A$ by sets $U_i \subseteq X$ with diameters $r_i < \epsilon$. The Hausdorff $\alpha$-measure of a set $A$ is defined as 
$$
H_{\alpha}(A) = \lim_{\epsilon \to 0}H_{\alpha, \epsilon}(A) = \sup_{\epsilon > 0}H_{\alpha, \epsilon}(A).
$$
Note that the limit exists.
Let $$H_{\mathrm{dim}}(A) = \inf\{\alpha : H_{\alpha}(A) = 0\}.$$ This is also equivalent to $$H_{\mathrm{dim}}(A) = \sup\{\alpha : H_{\alpha}(A) = \infty\}.$$
\end{definition} 

\section{A Family of Non-Uniquely Ergodic Interval Exchange Transformations}\label{sec:3}

In this section, we will construct a family of examples of non-uniquely ergodic transformations on $n$ intervals, and prove a series of properties, including obtaining an estimate for the Hausdorff dimension of the measures for such transformations. We will consider a Rauzy path with cycles that may depend on a parameter. In \cite{Fk}, Fickenscher provided a construction of a Rauzy path for a permutation with pairwise inversion of elements, where the first element is swapped with the last, and all other elements are swapped with their neighbors. By default, we consider the case where $n$ is even.

\subsection{Construction of the non-uniquely ergodic example}

$$ \pi =
\begin{Bmatrix}
1&2&3&...&\frac{n}{2}&\frac{n}{2}+1&...&n-2&n-1&n\\
n&3&2&...&\frac{n}{2}+1&\frac{n}{2}&...&n-1&n-2&1
\end{Bmatrix}
$$
\\
Let $a,c > 0$ and define the following Rauzy paths: 

$$\dot{\gamma}_{1, a} = 1^{n-2}0^a10^2,$$
$$\dot{\gamma}_{2, a} = 1^{n-3}0^a10^2,$$
$$...$$
$$\dot{\gamma}_{n-2, a} = 1^10^a10^2,$$
$$\gamma_{a, c} = 0\dot{\gamma}_{n-2, a}...\dot{\gamma}_{2, a}1^{c(n-1)}.$$
$$\\$$

The scheme of the Rauzy graph is presented as follows:

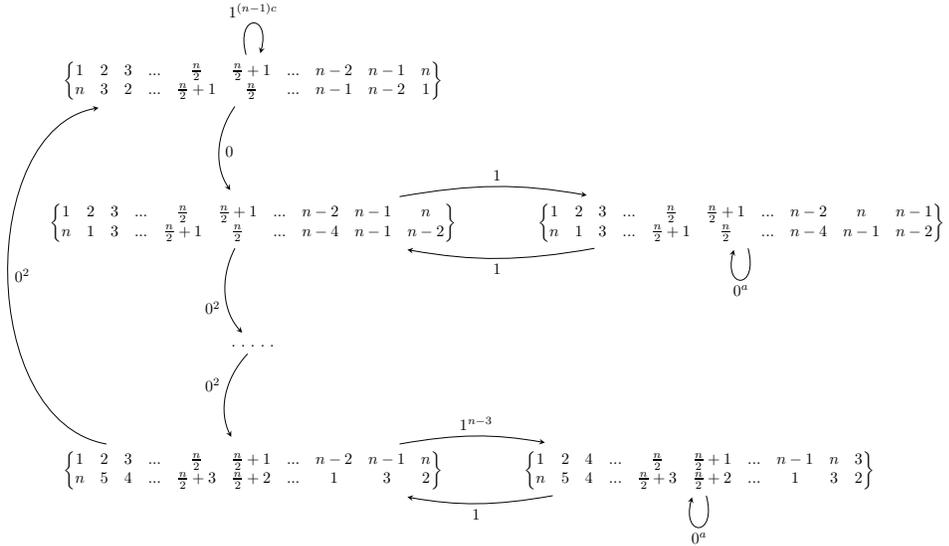
\scalebox{0.6}{
\begin{tikzpicture}[
    ->,
    >=stealth,
    shorten >=5pt,
    auto,
    node distance=2cm,
    semithick
]

\tikzset{
    state/.style={fill=white, draw=none, text=black, rounded corners, inner sep=4pt}
}

\node[state] (A) {$\begin{Bmatrix}
1&2&3&...&\frac{n}{2}&\frac{n}{2}+1&...&n-2&n-1&n\\
n&3&2&...&\frac{n}{2}+1&\frac{n}{2}&...&n-1&n-2&1
\end{Bmatrix}$};

\node[state] (B) [below=2cm of A] {$\begin{Bmatrix}
1&2&3&...&\frac{n}{2}&\frac{n}{2}+1&...&n-2&n-1&n\\
n&1&3&...&\frac{n}{2}+1&\frac{n}{2}&...&n-4&n-1&n-2
\end{Bmatrix}$};

\node[state] (C) [right=1.5cm of B] {$\begin{Bmatrix}
1&2&3&...&\frac{n}{2}&\frac{n}{2}+1&...&n-2&n&n-1\\
n&1&3&...&\frac{n}{2}+1&\frac{n}{2}&...&n-4&n-1&n-2
\end{Bmatrix}$};

\node[state] (D) [below=of B] {.\ .\ .\ .\ .};

\node[state] (F) [below=of D] {$\begin{Bmatrix}
1&2&3&...&\frac{n}{2}&\frac{n}{2}+1&...&n-2&n-1&n\\
n&5&4&...&\frac{n}{2}+3&\frac{n}{2}+2&...&1&3&2
\end{Bmatrix}$};

\node[state] (G) [right=1.5cm of F] {$\begin{Bmatrix}
1&2&4&...&\frac{n}{2}&\frac{n}{2}+1&...&n-1&n&3\\
n&5&4&...&\frac{n}{2}+3&\frac{n}{2}+2&...&1&3&2
\end{Bmatrix}$};

\path[->] (A) edge [bend right=35] node {0} (B);
\path[->] (A) edge [loop above] node {$1^{(n-1)c}$} (A);
\path[->] (B) edge [bend left=10] node {$1$} (C);
\path[->] (C) edge [bend left=10] node {$1$} (B);
\path[->] (C) edge [loop below] node {$0^a$} (C);
\path[->] (B) edge [bend right=35] node[swap] {$0^2$} (D);
\path[->] (D) edge [bend right=35] node[swap] {$0^2$} (F);
\path[->] (F) edge [bend left=10] node {$1^{n-3}$} (G);
\path[->] (G) edge [bend left=10] node {$1$} (F);
\path[->] (G) edge [loop below] node {$0^a$} (G);
\path[->] (F) edge [bend left=80] node[swap] {$0^2$} (A);

\end{tikzpicture}
}

Now let's construct the transition matrix corresponding to the cycle on the presented Rauzy diagram. Such a matrix represents the transition between the new and old lengths. The dependence of this matrix on the parameter arises from the fact that some nodes in the Rauzy path also depend on this parameter.

\begin{align*}
    \Theta_{\gamma_{a,c}}
        &= \Theta_{0}\Theta_{\dot{\gamma}_{n-2,a}} \cdots \Theta_{\dot{\gamma}_{2,a}}\Theta_{1^{c(n-1)}} \\
        &= \begin{pmatrix}
            1      & c     & c     &        &        & \cdots & c     \\
            0      & 2     & 2     & 1      &        & \cdots & 1     \\
            0      & a     & a+1   & 0      &        & \cdots & 0     \\
            0      & 1     & 1     & 2      & 2      & \cdots & 1     \\
            0      & 0     & 0     & a      & a+1    & \cdots & 0     \\
            \vdots &       &       & \ddots &        & \ddots & \vdots \\
            0      & \cdots& 1     & 1      & 2      & 2      & 1     \\
            0      & \cdots& 0     & 0      & a      & a+1    & 0     \\
            1      & c+1   & c+1   & \cdots & \cdots & \cdots & c+1
           \end{pmatrix}.
\end{align*}

If $n$ is odd, the first step in the path should be replaced by "$10$". In this case, the final matrix will have exactly the same form, but the $(n-1)$-th column will be equal to $\Theta_{\gamma_{a,c}}e_n + e_{n-1}$.

\subsection{Proof of the estimates for measures}

This section provides estimates for the measures of each interval of the transformation. The lemmas are technical results that offer insight into the frequency of visits to the subintervals, corresponding to the vector of lengths obtained from the IET constructed through a sequence of transition matrices.

\begin{lemm}
The following estimates hold for any $k \in \{3 , 5, \cdots, 2\floordown{\frac{n}{2}} - 1\}$: \\
1) $\frac{\lambda_k(I^j_k)}{\lambda_k(I^j)} > \frac{1}{2}$ \\
2) $\frac{\lambda_k(\cup_{i=2, i \neq k}^{n-1} I^j_i)}{\lambda_k(I^j)} < \frac{n}{2c^j}$
\end{lemm}

\begin{proof}
Base case: $$\Theta e_k = (c, 1, 0, \cdots, 0, 2, a+1, 1, 0, \cdots, 0, c+1)^T.$$
Let's check $x_k > \frac{1}{2}$:
\[
x_k = \frac{a+1}{2c+a+\frac{n}{2}+2} = \frac{pc+1}{(2+p)c+\frac{n}{2}+2} > \frac{1}{2}
\]
It is sufficient to show that
\[
2(pc+1) > (p+2)c+\frac{n}{2}+2
\]
\[
(p-2)c > \frac{n}{2}
\]
The last is true since $c > p > n$.

Let's check $\delta = x_2 + x_3 + \cdots + x_{k-1} + x_{k+1} + \cdots + x_{n-1} \le \frac{n}{2}c$:
\[
\delta = \frac{\frac{n}{2}}{(2+p)c+\frac{n}{2}+2} = \frac{n}{2(2+p)c+n+4} \le \frac{n}{2}c.
\]
This holds for any $p > 0$.

Inductive step: Assume it holds for $j+1$, let's show it for $j$.
We need to show that
$$2(\Theta x^{j+1})_k > |\Theta x^{j+1}|$$
Calculating the left-hand side:
$$
2(\Theta x^{j+1})_k = 2(a^{j}x^{j+1}_{k-1}+(a^{j}+1)x^{j+1}_k)
$$
Calculating the right-hand side:
\begin{multline*}
|\Theta x^{j+1}| = \\
2x^{j+1}_1 + (2c^{j}+\frac{n}{2}+a^{j}+1)\sum_{i=1}^{\ceilup{\frac{n}{2}} - 1}x^{j+1}_{2i}+(2c+\frac{n}{2}+a^{j}+2)\sum_{i=1}^{\ceilup{\frac{n}{2}} - 1}x^{j+1}_{2i+1}+(2c^{j}+\frac{n}{2})x^{j+1}_n
\end{multline*}

Let's estimate the expression $2(\Theta x^{j+1})_k - |\Theta x^{j+1}|$ from below:
\begin{multline*}
-2x^{j+1}_1 + \left((p-2)c^{j}-\frac{n}{2}-1\right)x^{j+1}_{k-1} + \left((p-2)c^{j}-\frac{n}{2}\right)x^{j+1}_{k} \\
- \left((p+2)c^{j}+\frac{n}{2}+1\right)\sum_{\substack{i=1 \\ i \neq k-1}}^{\lceil \frac{n}{2} \rceil - 1} x^{j+1}_{2i} \\
- \left((p+2)c^{j}+\frac{n}{2}+2\right)\sum_{\substack{i=1 \\ i \neq k}}^{\lceil \frac{n}{2} \rceil - 1} x^{j+1}_{2i+1} - \left(2c^{j}+\frac{n}{2}\right)x^{j+1}_n
\end{multline*}
We find a lower bound for this expression, which is achieved for the following distribution of values:
$$
x^{j+1}_k = \frac{1}{2}, x^{j+1}_3 = \delta^{j+1}, x^{j+1}_1 = x^{j+1}_2 = x^{j+1}_4 = \cdots = x^{j+1}_{n-1} = 0, x^{j+1}_n = \frac{1}{2} - \delta^{j+1},
$$
where $\delta^{j+1} = \sum_{i=2, i \neq k}^{n-1} x^{j+1}_i$ is the total weight of the intervals from the condition.
Thus, the lower bound is:
$$
\frac{1}{2}((p-2)c^{j}-\frac{n}{2})-((p+2)c^{j}+\frac{n}{2}+2)\delta^{j+1}-(2c^{j}+\frac{n}{2})(\frac{1}{2}-\delta^{j+1})=(\frac{p-4}{2})c^{j}-\frac{n}{2}-\frac{3pc^{j}+6}{p^2c^{j}},
$$
so it is sufficient to show that $(\frac{p-4}{2})c^{j}>\frac{n}{2}+\frac{3pc^{j}+6}{p^2c^{j}}$, which is true for large $p$.

Now let's prove the estimate for $\delta^{j}$.
Inductive step: Assume it holds for $j+1$, let's show it for $j$.
Note the following:
$$
|\Theta x^{j+1}|>(2c^{j}+\frac{n}{2}+a^{j}+2)x^{j+1}_k > \frac{1}{2}((p+2)c^{j}+\frac{n}{2}+2),
$$
Let's estimate $\delta^{j}$ from above:
$$
\delta^{j}<\frac{n-1+(n+pc^{j}+1)\frac{n}{2p^2c^{j}}}{\frac{1}{2}((p+2)c^{j}+\frac{n}{2}+2)}<\frac{n-1+\frac{n}{p}}{\frac{1}{2}(p+2)c^{j}} \le
 \frac{n}{2c^{j}},$$
 which is true for large $p$.
\end{proof}

\begin{lemm}
The following estimate is true: \\
1) $\frac{\lambda_1(I^j_1 \cup I^j_n)}{\lambda_1(I^j)} > \frac{1}{2}$ \\
2) $\frac{\lambda_1(\cup_{i=2}^{n-1} I^j_i)}{\lambda_1(I^j)} < \frac{n}{c^j}$
\end{lemm}

\begin{proof}
The base case is obvious, since:
$$\Theta e_1 = (1, 0, \cdots, 0, 1)^T.$$
Inductive step: Assume it holds for $j+1$, let's show it for $j$.
We need to show that
$$2((\Theta x^{j+1})_1+(\Theta x^{j+1})_n) > |\Theta x^{j+1}|$$
Calculating the left-hand side:
$$
(\Theta x^{j+1})_1+(\Theta x^{j+1})_n=2+(2c^{j}-1)\sum_{i=2}^{n-1}x^{j+1}_i=2c^{j}+1-(2c^{j}-1)x^{j+1}_1
$$
Calculating the right-hand side:

\begin{multline*}
|\Theta x^{j+1}| = \\
2x^{j+1}_1 + (2c^{j}+\frac{n}{2}+a^{j}+1)\sum_{i=1}^{\ceilup{\frac{n}{2}} - 1}x^{j+1}_{2i}+(2c^{j}+\frac{n}{2}+a^{j}+2)\sum_{i=1}^{\ceilup{\frac{n}{2}} - 1}x^{j+1}_{2i+1}+(2c^{j}+\frac{n}{2})x^{j+1}_n
\end{multline*}

Now let's compute the difference $2((\Theta x^{j+1})_1+(\Theta x^{j+1})_n)-|\Theta x^{j+1}|$:
\begin{multline*}
(4c^{j}+2) - (4c^{j}-4)x^{j+1}_1 - \left(2c^{j}+\frac{n}{2}+a^{j}+1\right)\sum_{i=1}^{\lceil \frac{n}{2} \rceil - 1} x^{j+1}_{2i} \\
- \left(2c^{j}+\frac{n}{2}+a^{j}+2\right)\sum_{i=1}^{\lceil \frac{n}{2} \rceil - 1} x^{j+1}_{2i+1} - \left(2c^{j}+\frac{n}{2}\right)x^{j+1}_n \\
> (4c^{j}+2) - (4c^{j}-4)x^{j+1}_1 - \left(2c+\frac{n}{2}\right)x^{j+1}_n \\
- \left(2c^{j}+\frac{n}{2}+a^{j}+2\right)\sum_{i=2}^{n-1} x^{j+1}_{i}
\end{multline*}
Let's estimate this expression from below by considering the worst case:
$$
(4c^{j}+2)-(4c^{j}-4)-\frac{(2c^{j}+\frac{n}{2}+a^{j}+2)n}{c^{j+1}}=6-\frac{((p+2)c^{j}+(n+4)/2)n}{p^2c^{j}}>0.
$$
Now let's prove the estimate for $\delta^{j}$.
Inductive step: Assume it holds for $j+1$, let's show it for $j$.
Note the following:
$$
|\Theta x^{j+1}|>(2c^{j}+\frac{n}{2})x^{j+1}_n>\frac{1}{2}c^{j}+\frac{n}{8},
$$
since $x_1+x_n > \frac{1}{2}$ and $x_n>x_1$, then $x_n > \frac{1}{4}$.
\begin{multline*}
\delta^{j} < \frac{(a^{j}+\frac{n}{2}+1)\sum_{i=2}^{n-1}x^{j+1}_i + (\frac{n}{2}-1)x^{j+1}_n}
               {\frac{1}{2}c^{j}+\frac{n}{8}} \\
           < \frac{(a^{j}+\frac{n}{2}+1)(n-2)}
               {p^2c^{j}\left(\frac{1}{2}c^{j}+\frac{n}{8}\right)}
               + \frac{n-2}{c^{j}+\frac{n}{4}}
           < \frac{1}{c^j} + \frac{n-2}{c^{j}} < \frac{n}{c^{j}}.
\end{multline*}
\end{proof}

The next theorem is a reformulated version of the result described in \\ Fickenscher's paper \cite{Fk}.
\begin{theorem}
Let $\lambda_j$ for $1 \le j \le n$ be the measures defined by the limiting distributions
\[
    \lim_{m \to \infty} \overline{\Theta}_{c_1,a_1}\overline{\Theta}_{c_2,a_2}\cdots\overline{\Theta}_{c_m,a_m} e_j.
\]
If $n$ is even, the following equalities hold:
\[
    \lambda_1 = \lambda_n, \quad \lambda_2=\lambda_3, \quad \dots, \quad \lambda_{n-2}=\lambda_{n-1}.
\]
If $n$ is odd, the equalities are:
\[
    \lambda_1 = \lambda_{n-1} = \lambda_n, \quad \lambda_2=\lambda_3, \quad \dots, \quad \lambda_{n-3}=\lambda_{n-2}.
\]
\end{theorem}
\begin{proof}
Let's prove for the even case. The proof for the odd case is similar. \\
The proof of the theorem follows from the combinatorial properties of the transition matrix $\Theta_{c,a}$. Note that 
$$ 
\Theta_{c, a}e_n = c\Theta_{c, a}e_1 + \Theta_{c, a}e_2 + \Theta_{c, a}e_3 + \cdots + \Theta_{c, a}e_n,
$$
$$
\Theta_{c, a}e_3 = \Theta_{c, a}e_2 + e_3,
$$
$$
\Theta_{c, a}e_5 = \Theta_{c, a}e_4 + e_5,
$$
$$
\cdots
$$
$$
\Theta_{c, a}e_{n-1} = \Theta_{c, a}e_{n-2} + e_{n-1}.
$$

Assume we have some sequence $(c_m,a_m)_{m=1}^{\infty}$, where $c_m > 0, a_m > 0$.
Consider the limit $\lambda_j =  \lim_{m \to +\infty} \overline{\Theta}_{c_1,a_1}\overline{\Theta}_{c_2,a_2}...\overline{\Theta}_{c_m,a_m} e_j$.
From the relations above, we get that $\lambda_1 = \lambda_n, \lambda_2=\lambda_3, \dots, \lambda_{n-2}=\lambda_{n-1}.$
This is true because the limit vector will be the same for each pair.
Since $M(T)$ is isomorphic to $\Lambda(\gamma)$, the transformation $T$ has no more than $\floordown{\frac{n}{2}}$ distinct invariant ergodic measures.
\end{proof}

\subsection{Construction of the non-uniquely ergodic family}
This section will provide the proof of Theorem \ref{onto}, which uses the estimates described in the lemmas. This theorem will allow us to construct $\floordown{\frac{n}{2}}$ distinct invariant ergodic measures.
\begin{proof}[Proof of Theorem \ref{onto}]
The proof of the system's non-unique ergodicity is based on demonstrating the existence of more than one invariant ergodic measure. We will construct $\floordown{\frac{n}{2}}$ such measures, $\lambda_3, \ldots, \lambda_{2\floordown{\frac{n}{2}}-1}, \lambda_n$, and prove that they are pairwise distinct.

To the measure $\lambda_j$ (where $j = 3, ..., 2\floordown{\frac{n}{2}}-1, n$) corresponds the vector of partition interval lengths given by the limit
$$\lim_{m \to +\infty} \overline{\Theta}_{c_1,a_1}\overline{\Theta}_{c_2,a_2}...\overline{\Theta}_{c_m,a_m} \mathbf{e}_j,$$
where $\mathbf{e}_j$ is the $j$-th standard basis vector.

All measures for this example are pairwise distinct.\\ If we consider $\alpha = (\alpha_{1}, \alpha_{2}, \dots, \alpha_{n})$ - the partition of the initial interval, then under the given constraints on $c_j$ and $a_j$ (namely, the condition $c^{i+1} = p a^i = p^2 c^i$ for a sufficiently large $p \geq n + 1$), we get that the relative frequency of the trajectory visiting the intervals corresponding to odd indices for the measure $\lambda_j$ will be predominantly concentrated on the interval with index $j$. More specifically, for each measure $\lambda_j$, the sum of the measures of intervals with odd indices will be greater than $\frac{1}{2}$, but the distribution of this mass among specific odd intervals will be different for different $j$. Since each measure $\lambda_j$ assigns the greatest weight to its "own" interval $\alpha_j$, the measures cannot coincide. It follows that $T$ is a non-uniquely ergodic transformation, and the constructed $\floordown{\frac{n}{2}}$ measures are pairwise distinct invariant ergodic measures.
\end{proof}

\section{Hausdorff Dimension of Measures}\label{sec:4}

\subsection{Adaptation of Frostman's Lemma}

We will use Frostman's lemma \cite{Fr} for the lower bound of the Hausdorff dimension.
\begin{lemm}[Frostman's Lemma]
 Let $B \subset [0, 1)$ be a Borel set. $H^s(B) > 0$ if and only if there exists a finite Radon measure $v$ on $B$ such that for any $x$ and $r>0$, we have $v(B(x, r)) \leq r^s$.
\end{lemm}

\begin{corollary}
If $\mu$ is a measure defined on $[0, 1)$ and $\epsilon_1, \epsilon_2, \epsilon_3, ...$ is a positive sequence tending to 0, such that ${{\epsilon_i}\over{\epsilon_{i+1}}} < C$ for some $C$, and for any $i$ we have $\mu(B(x, \epsilon_i)) \leq C(\epsilon_i)^{\alpha}$, it follows that $\dim_H(\mu) \geq \alpha$.
\end{corollary}

The adaptation of Frostman's Lemma to our specific example is formulated as follows:

\begin{lemm}\label{fl}
Let $\lambda_a$ and $\lambda_b$ be two distinct invariant ergodic measures for the map $T$.
If there exists a constant $C > 0$ and $\alpha \ge 0$ such that for any induction step $k$ and any $i \in \{1, \dots, n\}$ the inequality holds:
$$ \lambda_a(I_i^{(k)}) \le C \cdot \lambda_b(I_i^{(k)})^{\alpha},$$
then the Hausdorff dimension of measure $\lambda_a$ with respect to the metric $d_{\lambda_b}$ is bounded below:
$$\dim_H(\lambda_a, d_{\lambda_b}) \ge \alpha.$$
\end{lemm}

\begin{proof}
The proof follows directly from the application of the generalized Frostman's lemma (in the form of a corollary, as in J. Chaika's work) to our case, as it directly follows from the fact that $I_a^{(k)}$ and $I_b^{(k)}$ consist of repeating images and that any arbitrary interval can be effectively covered by a finite number of such intervals.
We choose as a sequence of scales the lengths of these intervals in the metric $\lambda_b$: $\epsilon_{i,k} = \lambda_b(I_i^{(k)}).$
Now we need to check the main condition of Frostman's lemma for our measure $\lambda_a$ and the chosen $\epsilon_{i,k}$. We need to show that for some constant $C'$:
    $$ \lambda_a(I_i^{(k)}) \le C' \cdot (\epsilon_{i,k})^\alpha $$
    Substituting our $\epsilon_{i,k}$:
    $$ \lambda_a(I_i^{(k)}) \le C' \cdot \lambda_b(I_i^{(k)})^\alpha $$
To find the best possible value of $\alpha$ satisfying the above condition for all $i$ and $k$, we must solve the inequality for $\alpha$:
    $$ C \cdot \lambda_b(I_i^{(k)})^{\alpha} \ge \lambda_a(I_i^{(k)}) $$
    $$ \lambda_b(I_i^{(k)})^{\alpha} \ge \frac{\lambda_a(I_i^{(k)})}{C} $$
    We obtain:
    $$ \alpha \le \log_{\lambda_b(I_i^{(k)})} \left( \frac{\lambda_a(I_i^{(k)})}{C} \right) = \log_{\lambda_b(I_i^{(k)})} (\lambda_a(I_i^{(k)})) - \log_{\lambda_b(I_i^{(k)})}(C). $$
    For this inequality to hold for all $i$ and $k$, $\alpha$ must be less than or equal to the infimum of the right-hand side over all $i,k$.
    We get that the lower bound for the dimension is defined as:
    $$\dim_H(\lambda_a, d_{\lambda_b}) \ge \sup(\alpha)$$
Thus, the lemma is proven.
\end{proof}

\subsection{Additional Estimates for Measures}

\begin{lemm}\label{l1}
The following relations are true: \\
1) $\lambda_1(I_3)=\lambda_1(I_5)=\cdots=\lambda_1(I_{2\floordown{\frac{n}{2}}-1})$ \\
2) $\lambda_1(I_2)=\lambda_1(I_4)=\cdots=\lambda_1(I_{2\floordown{\frac{n}{2}}})$ \\
3) $\forall k \in \{3,5,\cdots,2\floordown{\frac{n}{2}}-1\}:$
$$
\lambda_k(I_3)=\lambda_k(I_5)=\cdots=\lambda_k(I_{k-2})=\lambda_k(I_{k+2})=\cdots=\lambda_k(I_{2\floordown{\frac{n}{2}}-1}),
$$
$$
\lambda_k(I_2)=\lambda_k(I_4)=\cdots=\lambda_k(I_{k-3})=\lambda_k(I_{k+1})=\cdots=\lambda_k(I_{2\floordown{\frac{n}{2}}-2}).
$$
\end{lemm}

\begin{proof}
We first prove 1) and 2) together.
Base case - the vector before normalization is the first column of $\Theta$: $(1, 0, 0, 0, 0, 1)^T$.
Inductive step: Assume it holds for $j+1$, let's show it for $j$.
Let's compute all $x^j_2-x^j_4, \cdots, x^j_{2\floordown{\frac{n}{2}}-4}-x^j_{2\floordown{\frac{n}{2}}-2}$:
$$
x^j_{2k}-x^j_{2k+2} = \frac{x^{j+1}_{2k}+x^{j+1}_{2k+1}-x^{j+1}_{2k+2}-x^{j+1}_{2k+3}}{|\Theta|}
$$
Let's compute all $x^j_3-x^j_5, \cdots, x^j_{2\floordown{\frac{n}{2}}-3}-x^j_{2\floordown{\frac{n}{2}}-1}$:
$$
x^j_{2k+1}-x^j_{2k+3} = \frac{a^{j}(x^{j+1}_{2k}-x^{j+1}_{2k+2})+(a^{j}+1)(x^{j+1}_{2k+1}-x^{j+1}_{2k+3})}{|\Theta|}
$$
By the induction hypothesis, $x^{j+1}_{2k}-x^{j+1}_{2k+2}=x^{j+1}_{2k+1}-x^{j+1}_{2k+3}=0$, so:
$$
x^j_{2k}-x^j_{2k+2}=x^j_{2k+1}-x^j_{2k+3}=0.
$$
Let's prove 3):
Base case: it is sufficient to consider the vector:
$$\Theta e_k = (c, 1, 0, \cdots, 0, 2, a+1, 1, 0, \cdots, 0, c+1),$$
where $2, a+1$ correspond to positions $k-1, k$. Thus, the statement will be true, as even positions are equal to 1, and odd ones are 0.
Inductive step: Assume it holds for $j+1$, let's show it for $j$. \\
Let's compute all $x^j_2-x^j_4, \cdots, x^j_{k-3}-x^j_{k+1}, \cdots, x^j_{2\floordown{\frac{n}{2}}-4}-x^j_{2\floordown{\frac{n}{2}}-2}$:
$$
x^j_{2i}-x^j_{2i+2} = \frac{x^{j+1}_{2i}+x^{j+1}_{2i+1}-x^{j+1}_{2i+2}-x^{j+1}_{2i+3}}{|\Theta|}
$$
Let's compute all $x^j_3-x^j_5, \cdots, x^j_{k-2}-x^j_{k+2}, \cdots, x^j_{2\floordown{\frac{n}{2}}-3}-x^j_{2\floordown{\frac{n}{2}}-1}$:
$$
x^j_{2i+1}-x^j_{2i+3} = \frac{a^{j}(x^{j+1}_{2i}-x^{j+1}_{2i+2})+(a^{j}+1)(x^{j+1}_{2i+1}-x^{j+1}_{2i+3})}{|\Theta|}
$$
By the induction hypothesis, $x^{j+1}_{2i}-x^{j+1}_{2i+2}=x^{j+1}_{2i+1}-x^{j+1}_{2i+3}=0$, so:
$$
x^j_{2i}-x^j_{2i+2}=x^j_{2i+1}-x^j_{2i+3}=0.
$$
\end{proof}

\begin{lemm}
The following relations hold for $1 \le i \le \floordown{\frac{n}{2}}$: \\
1) $\lambda_1(I^k_{2i+1}) \le \frac{1}{p-1}\lambda_1(I^k_{2i}),$ \\ 2) $\lambda_1(I^k_{2i}) \le \frac{1}{c^k}\lambda_1(I^k_n),$ \\ 3) $\lambda_1(I^k_{2i}) \le \frac{1}{c^k-1}\lambda_1(I^k_1)$
\end{lemm}

\begin{proof}
Let's prove 3). The base case is obvious.
Inductive step: Assume it holds for $j+1$, let's show it for $j$.
It is sufficient to show it only for $\frac{x^j_1}{x^j_2}$, thanks to the lemma on equalities.
$$
\frac{x^j_1}{x^j_2} = \frac{c^{j} - (c^{j}-1)x^{j+1}_1}{\frac{n}{2}(x^{j+1}_2+x^{j+1}_3)+x^{j+1}_n}
$$
Let's show that the following is true:
$$
c^{j} - (c^{j}-1)x^{j+1}_1 \geq (c^{j}-1)(\frac{n}{2}(x^{j+1}_2+x^{j+1}_3)+x^{j+1}_n)
$$
Expanding and using the fact that $x^{j+1}_1+x^{j+1}_n = 1-(\frac{n}{2}-1)(x^{j+1}_2+x^{j+1}_3)$:
$$
1+\frac{1}{c^{j}-1} \geq 1+x^{j+1}_2+x^{j+1}_3,
$$
which is true since $x^{j+1}_2+x^{j+1}_3 \le \frac{1}{c^{j}} < \frac{1}{c^{j}-1}$.

Part 2) is proven analogously due to the similar structures of $I^k_1$ and $I^k_n$.

Let's prove part 1). The base case is obvious.
Inductive step: Assume it holds for $j+1$, let's show it for $j$:
$$
\frac{x^j_2}{x^j_3} = \frac{\frac{n}{2}(x^{j+1}_2+x^{j+1}_3)+x^{j+1}_n}{\frac{n}{2}(x^{j+1}_2+x^{j+1}_3)+x^{j+1}_n} \ge p-1
$$
Then we get
$$
x^{j+1}_n \ge ((p-1)a^j-\frac{n}{2})x^{j+1}_2 + ((p-1)(a^j+1)-\frac{n}{2})x^{j+1}_3
$$
Now it is sufficient to show that 
$$
\frac{x^{j+1}_n}{x^{j+1}_2} \ge ((p-1)a^j-\frac{n}{2}) + \frac{((p-1)(a^j+1)-\frac{n}{2})}{p-1}
$$
Which is equivalent to the inequality
$$
c^{j+1} = pa^j \ge ((p-1)a^j-\frac{n}{2}) + a^j + 1 - \frac{n}{2(p-1)}
$$
Then it is equivalent to
$$
\frac{n}{2} + \frac{n}{2(p-1)} \ge 1
$$
Which is obviously true.
\end{proof}

\begin{lemm}
The following relations hold for $i \in \{3, \cdots,2\floordown{\frac{n}{2}}-1\}$ and $j \neq \frac{i-1}{2}$:
\[
    \lambda_i(I^k_{2j}) > \lambda_i(I^k_{2j+1}).
\]
\end{lemm}

\begin{proof}
Prove this for $\lambda_3$. By Lemma~\ref{l1}, it suffices to show this for $j=2$.
First, we prove that
\[
    \lambda_3(I^k_{3}) > 2a^k\lambda_3(I^k_{4}),
\]
which holds by induction. The base case is obvious. Now we show the step:
\[
    x^j_4 \le \frac{(\frac{n}{2}+1)x^{j+1}_2+(\frac{n}{2}+1)x^{j+1}_3+x^{j+1}_n}{\lvert\Theta\rvert}.
\]
Note that
\[
    x^j_3 = \frac{a^jx^{j+1}_2+(a^j+1)x^{j+1}_3}{\lvert\Theta\rvert} > \frac{(\frac{n}{2}+1)a^{j-1}x^{j+1}_2+(\frac{n}{2}+1)a^{j-1}x^{j+1}_3+a^{j-1}x^{j+1}_n}{\lvert\Theta\rvert},
\]
since $p>n$ and $x^{j+1}_3 > x^{j+1}_n$.

Now, we will show that
\[
    \lambda_3(I^k_{4}) > \lambda_3(I^k_{5}).
\]
The base case is obvious.
Inductive step: Assume it holds for $j+1$, let's show it for $j$.
It is sufficient to show it only for the ratio $\frac{x^j_4}{x^j_5}$, thanks to the lemma on equalities.
\[
    \frac{x^j_4}{x^j_5} > 1,
\]
since $x^{j+1}_3 > 2a^jx^{j+1}_4$.
\end{proof}

\begin{lemm}
The order of decrease of subintervals for each measure is as follows:
$\lambda_1$: $I_n>I_1>I_2=I_4=\cdots=I_{2\floordown{\frac{n}{2}}-2}>I_3=I_5=\cdots=I_{2\floordown{\frac{n}{2}}-1}$ \\
$\lambda_3$: $I_3>I_n>I_1>I_2>I_4=\cdots=I_{2\floordown{\frac{n}{2}}-2}>I_5=I_7=\cdots=I_{2\floordown{\frac{n}{2}}-1}$ \\
$\lambda_5$: $I_5>I_n>I_1>I_4>I_2=\cdots=I_{2\floordown{\frac{n}{2}}-2}>I_3=I_7=\cdots=I_{2\floordown{\frac{n}{2}}-1}$ \\
$\cdots$ \\
$\lambda_{2\floordown{\frac{n}{2}}-1}$: \\
$I_{2\floordown{\frac{n}{2}}-1}>I_n>I_1>I_{2\floordown{\frac{n}{2}}-2}>I_2=\cdots=I_{2\floordown{\frac{n}{2}}-4}>I_3=I_5=\cdots=I_{2\floordown{\frac{n}{2}}-3}$
\end{lemm}

\begin{proof}
This follows from the previous lemmas.
\end{proof}

\begin{lemm} 
The following order will hold:
$$
b_{k,1} < b_{k,n} < b_{k,2} = b_{k,4} = \cdots = b_{k,2\floordown{\frac{n}{2}}} < b_{k,3} = b_{k,5} = \cdots = b_{k,2\floordown{\frac{n}{2}}-1}.
$$
\end{lemm}

\begin{proof}
The arrangement will depend on the initial order of the columns.
\end{proof}

\begin{lemm} 
The following estimate is true for $i \notin \{1, n\}$:
$$
\prod_{j=1}^k a^j < b_{k,i} < \prod_{j=1}^k 2a^j.
$$
\end{lemm}

\begin{proof}
To prove this, it is sufficient to consider $b_{k,2}$ and $b_{k,3}$, since $b_{k,2} < b_{k,3}$. We will show that
$$
\prod_{j=1}^k a^j < b_{k,2} < b_{k,3} < \prod_{j=1}^k 2a^j.
$$

This follows from the fact that
\begin{align*}
b_{k,3} &= c^kb_{k-1,1} + 2b_{k-1,2} + (a^k+1)b_{k-1,3} + (c^k+1)b_{k-1,n} \\
&\quad + \sum_j b_{k-1,j} \leq \left(2c^k + a^k + \frac{n}{2} + 2\right) b_{k-1,3} < 2a^kb_{k-1,3};
\end{align*}

\begin{align*}
b_{k,3} &= c^kb_{k-1,1} + 2b_{k-1,2} + (a^k+1)b_{k-1,3} + (c^k+1)b_{k-1,n} \\
&\quad + \sum_j b_{k-1,j} \geq (a^k+1)b_{k-1,3} > a^kb_{k-1,3};
\end{align*}
$$
b_{k,2} = c^kb_{k-1,1} + 2b_{k-1,2} + a^kb_{k-1,3} + (c^k+1)b_{k-1,n} + \sum_j b_{k-1,j} \geq a^kb_{k-1,3},
$$

from which we get

$$
\prod_{j=1}^k a^j < b_{k,2} < b_{k,3} < \prod_{j=1}^k 2a^j.
$$
\end{proof}

\begin{lemm} 
The following estimates are true:
$$
\prod_{j=1}^{k-1} c^j < b_{k,1} < \prod_{j=1}^k 2a^j; \\
\prod_{j=1}^{k} c^j < b_{k,n} < \prod_{j=1}^k 2a^j.
$$
\end{lemm}

\begin{proof}
This follows from the fact that
$$
b_{k,n} = c^kb_{k-1,1} + (c^k+1)b_{k-1,n} + \sum_j b_{k-1,j} > c^kb_{k-1,n};
$$
$$
b_{k,1} = b_{k-1,1} + b_{k-1,n} > b_{k-1,n};
$$
\end{proof}

\begin{lemm} 
The following inequality holds for all $i \in \{3,5,\cdots,2\floordown{\frac{n}{2}}-1,n\}$:
$$\lambda_i(O(I_i^{(k)})) > \frac{1}{n}$$
\end{lemm}

\begin{proof}
Since
$$
1 = \sum_{t=1}^n \lambda_i(O(I_t^{(k)})) = \sum_{t=1}^n b_{k,t} \lambda_i(I_t^{(k)}),
$$
then 
$$
\lambda_i(I^{(k)}) \sum_{t=1}^n b_{k,t} \frac{\lambda_i(I_t^{(k)})}{\lambda_i(I^{(k)})} = 1,
$$
from which we get
$$
\lambda_i(O(I_i^{(k)})) = \frac{b_{k,i} \frac{\lambda_i(I_i^{(k)})}{\lambda_i(I^{(k)})}}{\sum_{t=1}^n b_{k,t} \frac{\lambda_i(I_t^{(k)})}{\lambda_i(I^{(k)})}} = \frac{1}{1 + \frac{\sum_{t \neq i} b_{k,t} \frac{\lambda_i(I_t^{(k)})}{\lambda_i(I^{(k)})}}{b_{k,i} \frac{\lambda_i(I_i^{(k)})}{\lambda_i(I^{(k)})}}}
$$

Let's estimate the quantity from above:

$$
\frac{\sum_{t \neq i} b_{k,t} \frac{\lambda_i(I_t^{(k)})}{\lambda_i(I^{(k)})}}{b_{k,i} \frac{\lambda_i(I_i^{(k)})}{\lambda_i(I^{(k)})}} \leq \frac{b_{k,odd} \frac{n}{c^k}}{\frac{1}{4}b_{k,n}} + 2 \leq 4\frac{n}{p^{2k}} 2^{k}p^k + 2 < n - 1
$$

Therefore

$$
\lambda_i(O(I_i^{(k)})) > \frac{1}{1 + n - 1} = \frac{1}{n}
$$
\end{proof}

\subsection{Estimate of the Hausdorff Dimension}
In this section, we will prove several lemmas from which Theorem \ref{onto2} follows. This will allow us to conclude what the Hausdorff dimension of the measures obtained in Theorem \ref{onto} is.

\begin{lemm}\label{lem1}
Let $i$ and $j$ be distinct indices from the set $\{3,5,\cdots,2\floordown{\frac{n}{2}}-1,n\}$. Then the following estimate for the Hausdorff dimension holds:
\[ 
\dim_H(\lambda_i, d_{\lambda_j}) \geq \liminf_{k \to \infty} \log_{\lambda_j(I_i^{(k)})} (\lambda_i(I_i^{(k)})).
\]
\end{lemm}

\begin{proof}
An adaptation of Frostman's lemma \ref{fl} for interval exchange transformations provides the following general estimate:
$$
\dim_H(\lambda_i,d_{\lambda_j}) \geq \underset{1\leq t \leq n}{\min}\underset{k \to \infty}{\liminf} \log_{\lambda_j(I_t^{(k)})} (\lambda_iI_t^{(k)})).
$$
Now we need to understand for which $t$ the minimum is achieved. To do so, we analyze the asymptotic behavior of the logarithmic expression. Let us consider three cases.

Case 1: $i \neq n$ and $j \neq n$. From the preceding lemmas, it is known that for the measure $\lambda_i$, the measure of the interval $I_i^{(k)}$ is asymptotically the largest, whereas for the measure $\lambda_j$, the measure of the same interval $I_i^{(k)}$ is among the smallest. Therefore, the minimum is achieved when the argument is asymptotically maximal and the base is asymptotically minimal. Both of these conditions are met at $t=i$.

Case 2: $j=n$. The logic is analogous to the previous case. For the measure $\lambda_i$, the interval $I_i^{(k)}$ has the largest measure, while for the measure $\lambda_n$, its measure is among the smallest. Thus, the minimum is again achieved at $t=i$.

Case 3: $i=n$. We need to minimize the expression over $t$. \\ For $t \in \{3, 5, \dots, 2\lfloor\frac{n}{2}\rfloor - 1\}$ and any $j$ from the same set, the expressions $\log \lambda_j(I_t^{(k)})$ and $\log \lambda_n(I_t^{(k)})$ are asymptotically equivalent. However, it can be easy verified that $\lambda_n(I_n^{(k)}) > \lambda_j(I_n^{(k)})$. This implies that for $t=n$, the argument $\lambda_n(I_n^{(k)})$ of the logarithm is larger than for other $t$. Consequently, the minimum is also achieved at $t=i=n$. \\
In all cases, the minimum is achieved at $t=i$, which proves the lemma.
\end{proof}

\begin{lemm}\label{lem2}
For any distinct $i,j \in \{3,5,\cdots,2\floordown{\frac{n}{2}}-1,n\}$, the following estimate holds:
$$ 
\dim_H(\lambda_i,d_{\lambda_j}) \leq \underset{k \to \infty}{\liminf} \log_{\lambda_j(I_i^{(k)})} b^{-1}_{k,i}.
$$
\end{lemm}

\begin{proof}

First, we establish a connection between the dimension of the measure $\lambda_i$ and the dimension of the set $\LS O(I_i^{(k)})$. 
Recall that for any $k$, $\lambda_i(O(I_i^{(k)})) > \frac{1}{n}$. The set $\LS O(I_i^{(k)})$ is invariant under the map $T$, since the measure of the set $\lambda_i(I_i^{(k)}) \to 0$. The system is ergodic with respect to $\lambda_i$, then the invariant set $\LS O(I_i^{(k)})$ must have measure 0 or 1. Given the positive lower bound on the measure of each $O(I_i^{(k)})$, we conclude that $\LS O(I_i^{(k)})$ has full measure.

Now, let $s = \underset{k \to \infty}{\liminf} \log_{\lambda_j(I_i^{(k)})} b_{k,i}^{-1}$. We will show that $\dim_H(\lambda_i, d_{\lambda_j}) \le s$.
For this, it is sufficient to show that for any $\epsilon > 0$, we have $\dim_H(\lambda_i, d_{\lambda_j}) < s+\epsilon$. \\
By definition, there exists a subsequence of indices $\{k_m\}_{m=1}^\infty$ such that \\ $\lim_{m \to \infty} k_m = \infty$ and
$$ \lim_{m \to \infty} \log_{\lambda_j(I_i^{(k_m)})} b_{k_m,i}^{-1} = s $$
From this it follows that for any $\epsilon > 0$, we can choose this subsequence such that for all $m=1, 2, \dots$ we have:
$$ \log_{\lambda_j(I_i^{(k_m)})} b_{k_m,i}^{-1} < s + \epsilon $$
Now we show that $\dim_H(\LS O(I_i^{(k)}),d_{\lambda_j}) < s + \epsilon$:
We cover $ H^{s+\epsilon}_{\delta}(\LSM O(I_i^{(k_m)}))$ by covering each $O(I_i^{(k_m)})$ with a number $b_{k_m,i}$ of images of $I_i^{(k_m)}$. \\
Let's now prove that the series $\sum_{m=1}^\infty b_{k_m,i} (\lambda_j(I_i^{(k_m)}))^{s+2\epsilon}$ converges: \\
Since
$$ \log_{\lambda_j(I_i^{(k_m)})} b_{k_m,i}^{-1} < s + \epsilon $$
$$ \frac{-\ln(b_{k_m,i})}{\ln(\lambda_j(I_i^{(k_m)}))} < s + \epsilon $$
Then we get that 
$$ b_{k_m,i} < \frac{1}{(\lambda_j(I_i^{(k_m)}))^{s+\epsilon}} $$
$$ b_{k_m,i} \cdot (\lambda_j(I_i^{(k_m)}))^{s+\epsilon} < 1 $$
Thus, our series is majorized by the series $\sum_{m=1}^\infty (\lambda_j(I_i^{(k_m)}))^{\epsilon}$. This series converges.
\end{proof}

\begin{proof}[Proof of Theorem \ref{onto2}]
This follows directly from the proofs of Lemma \ref{lem1} and Lemma \ref{lem2}.
\end{proof}

After obtaining a generalization of J. Chaika's result for two measures from M. Keane's example, we can observe that the sequence of parameters defined in Theorem \ref{onto} also determines the Hausdorff dimension of the supports of the measures. This allows us to explicitly compute the Hausdorff dimension of the measures for a given sequence of parameters.

Next, we will provide a few more lemmas with estimates for our measures.


\begin{lemm} 
The following estimate holds:
$$
\dim_H(\lambda_i,d_{\lambda_j}) = \underset{k \to \infty}{\liminf} \log_{\lambda_j(I_i^{(k)})} b^{-1}_{k,i}
$$
\end{lemm}



\begin{proof}
Recall that $\lambda_i(O(I_i^{(k)})) = b_{k,i} \lambda_i(I_i^{(k)})$. From the previous estimates, we have the inequality:
$$
\frac{1}{n} \leq \lambda_i(O(I_i^{(k)})) \leq 1,
$$
which implies
$$
\frac{1}{n}b_{k,i}^{-1} \leq \lambda_i(I_i^{(k)}) \leq b_{k,i}^{-1}.
$$
Taking the logarithm with base $\lambda_j(I_i^{(k)})$, we obtain:
$$
\log_{\lambda_j(I_i^{(k)})} b_{k,i}^{-1} \leq \log_{\lambda_j(I_i^{(k)})} \lambda_i(I_i^{(k)}) \leq \log_{\lambda_j(I_i^{(k)})} b_{k,i}^{-1} - \log_{\lambda_j(I_i^{(k)})} n.
$$
Consider the limit of the term $-\log_{\lambda_j(I_i^{(k)})} n$ as $k \to \infty$. Since $\lambda_j(I_i^{(k)}) \to 0$, it follows that $-\log_{\lambda_j(I_i^{(k)})} n \to 0$. By applying the limits to the inequality, we get:
$$
\underset{k \to \infty}{\liminf} \log_{\lambda_j(I_i^{(k)})} (\lambda_i(I_i^{(k)})) = \underset{k \to \infty}{\liminf} \log_{\lambda_j(I_i^{(k)})} b^{-1}_{k,i}.
$$
\end{proof}




\begin{lemm} 
For $i \neq j$, it holds that
$$
\lambda_j(I_i^k) \leq \frac{1}{b_{k,i}}
$$
\end{lemm}

\begin{proof}
Since
$$
\sum_{t=1}^n b_{k,t} \lambda_j(I_t^k) = 1,
$$
it follows that
$$
b_{k,i}\lambda_j(I_i^k) \leq 1 \implies \lambda_j(I_i^k) \leq \frac{1}{b_{k,i}}
$$
\end{proof}

These estimates allow us to prove the corollary:

\begin{proof}[Proof of Corollary \ref{onto3}]
Since 
$$
\dim_H(\lambda_i,d_{\lambda_j}) = \underset{k \to \infty}{\liminf} \log_{\lambda_j(I_i^{(k)})} b^{-1}_{k,i} = \underset{k \to \infty}{\liminf} \frac{\log b^{-1}_{k,i}}{\log \lambda_j(I_i^{(k)})} \geq \underset{k \to \infty}{\liminf} \frac{\log b^{-1}_{k,i}}{\log b^{-1}_{k,i}} \geq 1.
$$
Since the Hausdorff dimension cannot exceed 1 in this context, we have equality.
\end{proof}

\subsection{Conclusion}

In this paper, we have conducted a detailed quantitative analysis of a family of maximally non-uniquely ergodic IETs by Fickenscher's construction. By developing novel combinatorial estimates and using Keane's method, we provided a proof for the existence of $\floordown{\frac{n}{2}}$ distinct invariant measures. Subsequently, following the methodology of J. Chaika, we established explicit formulas for their pairwise Hausdorff dimensions.

Several of the combinatorial estimates presented are sharp and provide significant insight into the underlying structure of this IET family. We believe these results lay a solid foundation for investigating other dynamical properties.

Looking forward, two main avenues for future research emerge from this work. First, a natural next step is to compute the Hausdorff dimension of the parameter space that defines non-uniquely ergodic IETs with maximal number of measures. The analytical tools and estimates introduced herein are expected to be valuable in pursuing this question. Second, can a construction analogous to Fickenscher's be realized for Interval Exchange Transformations with Flips.

\end{document}